\documentclass{article}
\usepackage[utf8]{inputenc}
\usepackage[english]{babel}

\usepackage[backend=bibtex,citestyle=numeric]{biblatex}
\usepackage{amsmath, amsthm, amssymb, tikz,enumerate,enumitem}
\usepackage{amsfonts}
\usepackage{fullpage}
\usepackage[enableskew,vcentermath]{youngtab}

\DeclareMathOperator{\orb}{orb}
\newcommand\beq{\begin{equation}}
\newcommand\eeq{\end{equation}}
\newcommand\bce{\begin{center}}
\newcommand\ece{\end{center}}
\newcommand\bea{\begin{eqnarray}}
\newcommand\eea{\end{eqnarray}}
\newcommand\ba{\begin{array}}
\newcommand\ea{\end{array}}
\newcommand\ben{\begin{enumerate}}
\newcommand\een{\end{enumerate}}
\newcommand\bit{\begin{itemize}}
\newcommand\eit{\end{itemize}}
\newcommand\brr{\begin{array}}
\newcommand\err{\end{array}}
\newcommand\bt{\begin{tabular}}
\newcommand\et{\end{tabular}}

\newtheorem{theorem}{Theorem}[section]
\newtheorem*{theorem*}{Theorem}

\newtheorem{proposition}[theorem]{Proposition}
\newtheorem{lemma}[theorem]{Lemma}

\newtheorem{definition}[theorem]{Definition}
\newtheorem{example}[theorem]{Example}

\newtheorem{defn}[theorem]{Definition}

\newtheorem{remark}[theorem]{Remark}

\newtheorem{observation}[theorem]{Observation}

\newcommand{\todo}[1]{\vspace{5 mm}\par \noindent
\marginpar{\textsc{ToDo}} \framebox{\begin{minipage}[c]{0.95
\textwidth}
 \end{minipage}}\vspace{5 mm}\par}

\begin{document}
\title{The poset of King permutations on a cylinder}
\author{Eli Bagno, Estrella Eisenberg, Shulamit Reches and Moriah Sigron}


\maketitle

\begin{abstract}
A permutation $\sigma=[\sigma_1,\dots,\sigma_n] \in S_n$ is called a {\em cylindrical king permutation} if $ |\sigma_{i+1}-\sigma_{i}|>1$ for each  $1\leq i \leq n-1$ and  $|\sigma_1-\sigma_n|>1$. 
The name comes from the the way one can see these permutations as describing 
locations of $n$ kings on a chessboard of order $n\times n$ in such a way that (each row and each column contains exactly one king and) no two kings are attacking each other, with the additional condition that 
a king can move off a certain row and reappear at the beginning of that row. 

In a recent paper, we dealt with the more general set of 'king permutations' i.e. the ones which satisfy only the first of the two conditions above. This set constitutes a poest under the well known containment relation on permutations. 

In this article we investigate the sub-poset of the cylindrical king permutations and its structure.  
We examine those cylindrical king permutations whose downset is as large as possible in the upper ranks. We use a modification of Manhattan distance of the plot of a permutation and some of its applications to the cylindrical context to find a criterion for such a permutation to be $k-$ prolific. 
One of our main results is that the maximal gap between two permutations in the poset of cylindrical permutations is $4$. 

\end{abstract}

\section{Introduction}

A {\em king permutation} is  a permutation $\sigma=[\sigma_1,\dots,\sigma_n]\in S_n$ such that for each $1 \leq i  \leq n-1$, $|\sigma_{i+1}-\sigma_i|>1$. We denote the set of the king permutations of order $n$ by $K_n$. For example: 
$K_1=S_1, K_2=K_3=\emptyset$, $K_4=\{[3142],[2413]\}$. 
The name "king permutations" refers to the way one can see these permutations as describing locations of $n$ kings on a chessboard of order $n\times n$ in such a way that (each row and each column contains exactly one king and) no two kings are attacking each other. In \cite{BERSKINGS} we investigated the structure of the poset of king permutations with respect to containment.

In our recent paper \cite{BERS1}, we considered a celebrated variant of the chess game, namely the cylindrical chess, in which the right and left edges of the board are imagined to be joined in such a way that a piece can move off one edge and reappear at the opposite edge. Let us call this board a {\it cylindrical board}. 

In the language of permutations, non-attacking kings that may move on a $n \times n$ cylindrical board can be described by the following subset of $K_n$:

\begin{definition}
Let $CK_n$ be defined as

$$CK_n=\{\sigma=[\sigma_1,\dots,\sigma_n] \in S_n: |\sigma_i-\sigma_{i+1}|>1,|\sigma_1-\sigma_n|>1,1\leq i \leq n-1\}.$$ 
An element of $CK_n$ will be called a {\em cylindrical king permutation}.

\end{definition}
For example: $CK_1=S_1, CK_2=CK_3=CK_4=\emptyset$, and $$CK_5=\{[31425],[14253],[42531],[25314],[53142],[24135],[41352],[13524],[35241],[52413]\}.$$

For $\sigma\in S_n$ and $\pi \in S_k$, $(k<n)$, we say that $\sigma$ {\it contains} $\pi$ if there is a subsequence of (the one line notation of) $\sigma$ that is order-isomorphic to that of $\pi$. 
If this is the case, we write $\pi \prec \sigma$. In this study we are interested in the sub-poset $\mathbb{CK}=\bigcup\limits_{n\geq 1}CK_n$  of the poset $\mathbb{K}=\bigcup\limits_{n\geq 1}K_n$ with respect to the containment relation. 

 A {\it bond} in a permutation $\pi\in S_n$ is a length $2$ consecutive sub-sequence of adjacent numbers (i.e. interval of length $2$ in the terminology of \cite{bevan}). Note that a king permutation is actually a permutation without bonds.  This point of view can be used also to describe the set of cylindrical king permutations. 
To this aim we extended in \cite{BERS1} the definition of bonds to cyclic bonds.
\begin{defn}\label{definition of blocks}
Let $\sigma =[\sigma_1,\ldots,\sigma_n]\in S_n$ and let $i \in [n-1]$. We say that the pair $(\sigma_i,\sigma_{i+1})$ is a {\em (regular) bond} in $\sigma$ if $\sigma_{i}-\sigma_{i+1} =\pm 1$. If $\sigma_{n}-\sigma_1=\pm 1$ then we say that the pair $(\sigma_n,\sigma_1)$ is an {\em edge bond} of $\sigma$. In general, adopting the convention that $\sigma_{n+1}=\sigma_1$, we say that $(\sigma_i,\sigma_{i+1})$ is a cyclic bond if it is a regular or an edge bond.
\end{defn}

According to this new definition, a permutation is a cylindrical king permutation if and only if it has no
cyclic bonds.\\

\subsection{Overview and main results}
In this paper we investigate the poset  $\mathbb{CK}$ of the cylindrical king permutations under the containment relation and its structure.
 In addition, we examine those cylindrical king permutations whose downset is as large as possible in the upper ranks. 

Here is an example of the downset of cylindrical king permutations of the permutation: $[5246173]$\\

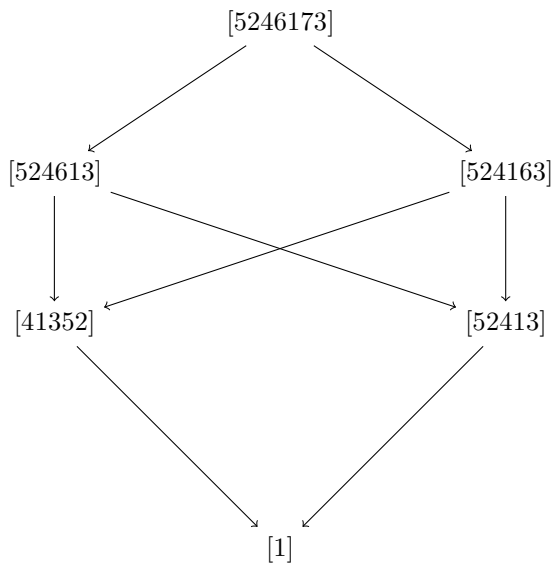
\begin{tikzpicture}
       \tikzstyle{every node} = [rectangle]
        \node (name) at (0,-1) {$$}; 
1

         \node (1) at (0,0) {$[1]$};
         \node at (0,0.5)  {};
        \node (41352) at (-3,3) {$[41352]$};
       \node at (-3.2,3.5)  {};
        
        \node (52413) at (3,3) {$[52413]$};
        \node at (3.2,3.5)  {};
      
        \node (524163) at (3,5) {$[524163]$};
        \node at (3,5.5)  {};
        
        \node (524613) at (-3,5) {$[524613]$};
        \node at (-3,5.5)  {};
        
        \node (5246173) at (0,7) {$[5246173]$};
        \node at (0,7.5) {} ;
        
        \foreach \from/\to in {41352/1,52413/1,524163/52413,524613/52413,524613/41352,524163/41352,5246173/524613,5246173/524163}
            \draw[->] (\from) -- (\to);
    \end{tikzpicture}

Recall that one can view a permutation via its graphical representation,
and let us identify a permutation $\sigma=[\sigma_1,\dots,\sigma_n]\in S_n$  with its plot, i.e. the set of all lattice points of the form $(i,\sigma_i)$ where $1 \leq i \leq n$. From this point of view, we can recognize king permutations, as well as cylindrical king permutations, by considering the distances between the points of their plots.
The {\it Manhattan} or {\it taxicab} distance is defined as follows:

\begin{definition}
Let $\sigma \in S_n$ and let $i,j\in [n]$. The (Manhattan) distance between the $i-th$ and the $j-th$ entries is defined to be the $L_1$ distance between the corresponding points in the plot of $\sigma$:

$$d_{\sigma}(i,j)=||(i,\sigma_i)-(j,\sigma_j)||_1=|i-j|+|\sigma_i-\sigma_j|.$$
\end{definition}
In \cite{BHT}, the concept of {\it breadth} of a permutation $\sigma \in S_n$ is defined to be the minimal distance between any two distinct entries of $\sigma$. Explicitly:
 $$br(\sigma)=min_{i,j\in [n],i\neq j}d_{\sigma}(i,j). $$

\begin{example}
For $\sigma=[5371426]$ we have $d_\sigma(2,5)=4$, $d_\sigma(1,4)=7$, $d_\sigma(1,2)=3$ and its breadth is $br(\sigma)=3$
\end{example}
Using the above definition we can characterize a king permutation as follows: 
\begin{observation}
 A permutation $\sigma $ of order $n>1$ is a king permutation if and only if its breadth is at least $3$. 
\end{observation}
In order to characterize cylindrical king permutations in a similar way, we modify in section \ref{two} the concepts of (Manhattan) distance and breadth to cyclic (Manhattan) distance and cyclic breadth  respectively.

For a permutation $\sigma \in S_n$ and $1 \leq k\leq n-1$, the maximal possible number of distinct permutations $\pi \in S_{n-k}$ which are contained in $\sigma$ is $\binom{n}{k}$.  
Indeed, such a permutation $\pi$, is obtained from $\sigma$ by deletion of $k$ entries from the one-line notation
of $\sigma$. Following \cite{BHT}, our interest is in those permutations $\sigma \in CK_n$ which contain maximally many patterns $\pi \in CK_{n-k}$.

 To this aim, 
we use and modify the following concept of {\it $k$-prolific permutation} which was defined 
 in \cite{BHT}:
 
\begin{definition}
A permutation $\sigma \in S_n$ is called $k$-prolific if $|\{\pi\in S_{n-k} | \pi\prec \sigma\}|=\binom{n}{k}. $ 
\end{definition}
In other words, a permutation $\sigma=[\sigma_1,\dots,\sigma_n]$ is  $k$-prolific if each $(n-k)-$subset of entries of $\sigma$ forms a unique pattern in $S_{n-k}$. It is known (see \cite{BERSKINGS}) that $\sigma \in S_n$ is $1-$prolific if and only if $\sigma$ is a king permutation, and this happens exactly when $br(\sigma)\geq 3$. 
We are interested in characterizing the cylindrical king permutations which are $k-$prolific in $CK_n$. I.e the cylindrical king permutations which contain a maximal number of distinct cylindrical king patterns.

 In Section \ref{two} we first prove that omitting a single entry from a permutation may decrease the cyclic breadth by at most one, and that for $\sigma \in CK_n$, there are $n$ distinct permutations $\sigma' \in CK_{n-1}$ such that $\sigma' \prec \sigma$ if and only if the cyclic breadth of $\sigma$ is greater than $3$ (See Proposition \ref{single entry} and Proposition \ref{cbr>3} respectively).  Then, using these propositions and our modifications to the concept of $k-$prolific, we get the two main results of this section, (which are generalizations of Theorem 2.25 in \cite{BHT} to king permutation and cylindrical king permutations respectively):
 
 \begin{enumerate}
     \item  A permutation $\sigma$ is $k$-prolific in $K_n$ if and only if  its {\bf breadth} is greater or equal to $k+3$. (See Theorem \ref{criterion for k prolific in Kn})
     
 \item A permutation $\sigma$ is $k$-prolific in $CK_n$ if and only if its {\bf cyclic breadth} is greater or equal to $k+3$. (See Theorem \ref{criterion for k prolific in CKn}).
\end{enumerate}
In section \ref{structure} we investigate the structure of the poset of cylindrical king permutations. This section is a direct continuation of our previous paper \cite{BERSKINGS}, which discusses the structure of the poset of king permutations. 
Since the cylindrical king permutations is a distinguished subset of the king permutations, we examine which  features of the poset of the cylindrical king  are inherited from the kings poset and which features are different.
One of our results is identifying the building blocks of the poset of cylindrical king permutations:
\begin{itemize}
    \item For every $\sigma  \in CK_n$ ($n \geq 5$) there is some $\pi \in CK_5$ such that $\pi \preceq \sigma$.  (See Theorem \ref{Each king contains 31425 or 24135}).
\end{itemize}
Our main result in this section is that in contrast to the case of king permutations in which the maximal gap between two permutations is 3, in $\mathbb{CK}$ the maximal gap is $4$:
\begin{itemize}
    \item Let $\sigma, \pi \in \mathbb{CK}$ be such that $\pi \prec \sigma$ and $|\sigma|-|\pi|> 5$.
Then there is some $\tau \in \mathbb{CK}$ such that $\pi \prec \tau \prec \sigma$ and $|\sigma|-|\tau|\leq 4$. (See Theorem \ref{antisematic}). 
\end{itemize}
Finally, in section \ref{direction} we present directions for further research in this area.

  \section{Cyclic breadth and $k$-prolific cylindrical king permutations} 
  \label{two}

We are interested in characterizing the cylindrical king permutations and analyzing some properties of their downsets in $\mathbb{CK}$. In particular, in this section, we focus on those cylindrical king permutations which contain as many distinct patterns as possible in the poset $\mathbb{CK}$. There is a close relationship between the breadth of a permutation and the number of distinct patterns it contains. Moreover, as mentioned above, one way to identify a king permutation is by verifying that its breadth is greater than or equal to $3$. In order to be able to identify cylindrical king permutations in a similar way, and especially those whose downsets in $\mathbb{CK}$
 have maximal size, we modify the definitions of the {\it (Manhattan) distance} and the {\it breadth}, introducing the {\it cyclic (Manhattan) distance} and the {\it cyclic breadth}
of a permutation $\sigma \in S_n$.

\begin{figure}[!ht]
    \centering
    \begin{tikzpicture}[cap=round,line width=2pt]
  \draw  (0,0) circle (2cm);

   \foreach \angle in
    {0, 45, 90, 135, 180, 225, 270, 315 }
  {
    \draw[line width=2pt] (\angle:1.9cm) -- (\angle:2.1cm);
  }


\node[right] at (2.1cm,0){ ${\bf\sigma_3=4}$}; 
\node[right] at (45:2.1cm){ $\sigma_2=6$}; 
\node[left] at (45:1.9cm){\tiny $i=2$};
\node[above] at (90:2.1cm){ ${\sigma_1=2}$}; 
\node[above left] at (135:2cm){ $\sigma_8=7$}; 
\node[left] at (180:2.1cm){ ${\bf\sigma_7=3}$}; 
\node[right] at (180:1.9cm){\tiny $j=7$};
\node[below left] at (225:2.1cm){ ${\bf\sigma_6=8}$}; 
\node[below] at (270:2.1cm){ ${\bf\sigma_5=5}$}; 
\node[below right] at (315:2.1cm){ ${\bf\sigma_4=1}$}; 

\end{tikzpicture}

\caption{$\sigma=[26415837]$  located on a circle in such a way that $\sigma_1$ occupies the north pole.
$cd_\sigma(2,7)=|6-3|+||7-2||=6$ where $||7-2||=min\{7-2,8+2-7\}$.}
    \label{S1}
\end{figure}
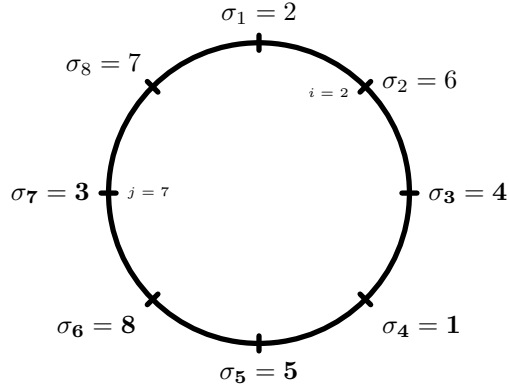

Along this section we consider a permutation $\sigma=[\sigma_1,\ldots,\sigma_n]$ as if it was located on a circle in such a way that $\sigma_1$ occupies the north pole (See figure \ref{S1}).

\begin{definition}
For $1 \leq i<j \leq n$, we denote by 
$((i,j))$ the shorter path on the circle leading from $i$ to $j$ which its length is the minimum of $j-i$ and $n-j+i$. 
We also set $||j-i||=|((i,j))|$, i.e. 
$||j-i||=min(j-i,n-j+i)$ 
\end{definition}
\begin{definition}
Let $\sigma \in S_n$ and let $i<j\in [n]$. Let the cyclic  (Manhattan) distance between the $i-th$ and the $j-th$ entries be defined as:

$$cd_{\sigma}(i,j)=|\sigma_j-\sigma_i|+||j-i||.$$

The {\em cyclic breadth} of $\sigma \in S_n$ is then defined to be:
$$cbr(\sigma)=min_{i,j\in [n],i\neq j}cd_{\sigma}(i,j).$$
\end{definition}

\begin{example}
\label{exa}
Let $\sigma = [724915836] \in S_9$. Then $cd_\sigma(1, 2) = (7-2)+||2-1|| = 6$, $cd_\sigma(2, 5) = (2-1)+||5-2|| = 4$ and $cd_\sigma(1, 9) =(7-6)+||9-1||=2$. The breath of $\sigma$ is $cbr(\sigma)=2$ 
\end{example}

We compare now between the regular (Manhattan) distance and the cyclic (Manhattan) distance as follows: 
\begin{observation}
\begin{itemize}
    \item For each $\sigma \in S_n$ and for each $i,j \in [n]$, $cd_{\sigma}(i,j) \leq d_{\sigma}(i,j)$
and so $cbr(\sigma) \leq br(\sigma)$.
\item For every $i,j \in [n]$, we have $d_\sigma(i,j)=d_{\sigma^{-1}}(\sigma_i,\sigma_j)$ and thus $br(\sigma)=br(\sigma^{-1})$ which means that the regular (Manhattan) breadth is closed under inverse.
\item The cyclic breadth $cbr(\sigma)$ is not always equal to $cbr(\sigma^{-1})$. 
For example, if $\sigma = [72415836] \in S_8$  then $\sigma^{-1}=[42735816]$, so that  $cbr(\sigma)=cd_{\sigma}(1,8)=2$ but  $cbr(\sigma^{-1})=cd_{\sigma^{-1}}(1,2)=3$.
\end{itemize}
\end{observation}

Using our new definitions, we can identify the cylindrical king permutation as follows: 
\begin{observation}
\label{obs1}
A permutation $\sigma $ of order $n>1$ is a cylindrical king permutation if and only if its cyclic breadth is greater than or equal to $3$.
\end{observation}
Omitting  a  single  entry  from  a  permutation  $\sigma \in S_n$ (and then standardising) produces a permutation $\sigma' \in S_{n-1}$ s.t. $\sigma' \prec \sigma$. 
In Proposition 2.24 of \cite{BHT}, the authors proved a proposition which partially controls the change in breadth of permutations in $S_n$: 
omitting a single entry from a permutation decreases the breadth by at most one.
We now present and prove a proposition that extends the above result, considering the cyclic breadth rather than the regular breadth. This enables us to characterize the cyclic breadth of permutations $\sigma'$ which are contained in a certain permutation $\sigma$.


\begin{proposition} 
\label{single entry}
Omitting a single entry from a permutation may decrease the cyclic breadth by at most one.
\end{proposition}

 \begin{proof}
Let $\sigma=[\sigma_1,\sigma_2,\ldots,\sigma_n] \in S_n$ and let $1\leq i<j\leq n$ be such that without loss of generality $\sigma_i<\sigma_j$. 
Recall that $cd_\sigma(i,j)$ is composed of the {\it position part} $||j-i||$ and the {\it value part} $|\sigma_j-\sigma_i|$. Let $\sigma'$ be the permutation obtained from $\sigma$ by omitting an element $\sigma_k$ for some $k \notin \{i,j\}$. We discuss the effect of omitting $\sigma_k$ on each one of the parts of $cd_\sigma(i,j)$. 

Observe that if $\sigma_k>\sigma_j$ or $\sigma_k<\sigma_i$ then this omission does not affect the value part of $cd_{\sigma}(i,j)$.  On the other hand, if $\sigma_i < \sigma_k < \sigma_j$ then, its omission decreases the value part of $cd_{\sigma}(i,j)$ by $1$.\\ 
Similarly, the omission of $\sigma_k$ affects the position part if and only if $k\in ((i,j))$. In case there is an effect, the difference is $1$. 

We conclude, that the omission of $\sigma_k$ decrease $cd_{\sigma}(i,j)$ of $\sigma$ by $2$ if and only if $\sigma_i<\sigma_k<\sigma_j$ and $k \in ((i,j))$, in all other cases the difference is $0$ or $1$.






Now, let $i<j$ such that $cbr(\sigma)=cd_{\sigma}(i,j)$. Assume to the contrary that there is $k \in [n]$ such that $cbr(\sigma')=cbr(\sigma)-2$, where $\sigma'$ is obtained from $\sigma$ by omitting $\sigma_k$. This means that $k\in ((i,j))$ and  $\sigma_i<\sigma_k<\sigma_j$.   We have now:
$$cbr(\sigma)\leq cd_{\sigma}(i,k)=\sigma_k-\sigma_i+||k-i||
< \sigma_j-\sigma_i+ ||j-i||=cd_{\sigma}(i,j),$$

a contradiction.

    
    
  We showed that if $cbr(\sigma)=cd_{\sigma}(i,j)$, then the omission of any element that is not $\sigma_i$ or $\sigma_j$ decreases the breadth by at most one.
    Finally, we show that for $cbr(\sigma)=cd_{\sigma}(i,j)$, an omission of one of the elements $\sigma_i$ or $\sigma_j$ also decreases the cyclic breadth by at most one. Without loss of generality let us consider the omission of $\sigma_j$.
    
    Indeed, for every $m,l \neq j$ we have $cbr(\sigma)=cd_{\sigma}(i,j) \leq cd_{\sigma}(m,l)$. If $cd_{\sigma}(i,j) < cd_{\sigma}(m,l)$ then the omission of $\sigma_j$ decreases $cd_{\sigma}(m,l)$ by at most $2$ to a value that is greater or equal to $cd_{\sigma}(i,j)-1$.  If $cd_{\sigma}(i,j) = cd_{\sigma}(m,l)$ then $cbr(\sigma)=cd_{\sigma}(m,l)$, and due to the previous discussion, the omission of $\sigma_j$ decreases $cd_{\sigma}(m,l)$ by at most one. Therefore if $\sigma'$ is obtained from $\sigma$ by omitting $\sigma_j$ then also $cbr(\sigma') \geq cbr(\sigma)-1$.

\end{proof}

\begin{remark}
Note that omitting an element $\sigma_k$ might also {\bf increase} the breadth of a permutation $\sigma$. For example if $\sigma=[351246]$ then $cbr(\sigma)=cd(3,4)=2$. If we omit $\sigma_3=1$, we get $\sigma'=[24135]$ with 
$cbr(\sigma')=3$.\\
\end{remark}


In \cite{Coleman}, Coleman demonstrates that maximising the distance between every two elements in a permutation tends to increase
the number of distinct sub-permutations (patterns).

In the proof of the next proposition and of several propositions and theorems in the next section, we  use the concept of cyclic separators and the notation $\nabla_i^*(\sigma)$.

Let $\nabla^*_i(\sigma)$ be the  permutation obtained by omitting the value $i$ from $\sigma$ and standardize. For example if $\sigma=[31425]$ then $\nabla^*_2(\sigma)=[2134]$ and $\nabla^*_3(\sigma)=[1324]$.

 The definition of cyclic separators appears in the Appendix of this paper.  

The next proposition characterizes the permutations $\sigma \in CK_n$ which contain $n$ different patterns $\sigma' \prec \sigma$ s.t. $\sigma' \in CK_{n-1}$.

\begin{proposition}\label{cbr>3}
 For $\sigma \in CK_n$, there are $n$ distinct permutations $\sigma' \in CK_{n-1}$ such that $\sigma' \prec \sigma$ if and only if  $cbr(\sigma)>3$.  
\end{proposition}
\begin{proof}
Let $\sigma \in CK_n$ and suppose that there are $n$ distinct permutations  $\sigma' \in CK_{n-1}$ such that $\sigma' \prec \sigma$, which means that by omitting any entry of $\sigma$, we get a cylindrical king permutation. Since $\sigma \in CK_n$, according to observation \ref{obs1} $cbr(\sigma)\geq 3$. Assume to the contrary that $cbr(\sigma)=3$, thus there are entries $i$ and $j$, $i<j$ such 
that $cbr(\sigma)=cd_{\sigma}(i,j)=3$, which means that $||j-i||=1$ and $|\sigma_j-\sigma_i|=2$ or vice versa. Assume without loss of generality that $\sigma_i<\sigma_j$.
\begin{itemize}
    
    \item If $||j-i||=1$ and $\sigma_j-\sigma_i=2$ then $\sigma_i+1$ is a horizontal cyclic separator, so that $\sigma'=\nabla_{\sigma_{i}+1}^*(\sigma)\notin CK_{n-1}$, a contradiction.

    
    \item If $||j-i||=2$ and $\sigma_j-\sigma_i=1$ then $\sigma_{i+1}$ is a vertical cyclic separator and again $\sigma'=\nabla_{\sigma_{i+1}}^*(\sigma)\notin CK_{n-1}$, a contradiction.
    
\end{itemize} 

Now, suppose that $cbr(\sigma)>3$. By corollary 7 of \cite{H2}, a permutation $\sigma \in S_n$ has maximum number of patterns $\sigma' \prec \sigma $ with  $\sigma' \in S_{n-1}$, if and only if it has no bonds. The permutation $\sigma \in CK_n$  has no bonds as a cylindrical king, thus there are $n$ distinct permutations $\sigma'\in S_{n-1}$ such that $\sigma' \prec \sigma$. We show now that each one of them is a cylindrical king permutation.  According to Proposition \ref{single entry}, omitting a single entry from a permutation decreases the cyclic breadth by at most one. Thus, each $\sigma'\in S_{n-1}$  such that $\sigma' \prec \sigma$  satisfies $cbr(\sigma')\geq 3$ and by observation \ref{obs1}, $\sigma'\in CK_{n-1}$ .  
\end{proof}

Recall that  a permutation $\sigma \in S_n$ is  $k$-prolific if each $(n-k)-$subset of entries of $\sigma$ forms a unique pattern. In other words, a permutation $\sigma \in S_n$ is  $k$-prolific  if it contains the maximal possible number of distinct permutations $\pi \in S_{n-k}$ such that $\pi\prec\sigma$.  
The breadth of a permutation affects its being a $k-$ prolific as shown in Theorem 2.25 in \cite{BHT}:  
$\sigma$ is $k$-prolific if and only if $br(\sigma)\geq k+2$. This idea can be extended to the poset of king permutations as well as to the poset of cylindrical king permutations. This is the context of the following definition. 
\begin{definition}
We say that a permutation $\sigma\in K_n, (CK_n)$ is {\it $k$-prolific in $K_n, (CK_n)$} if each $(n-k)$-subset of the entries of $\sigma=[\sigma_1,\dots,\sigma_n]$ forms a unique pattern in $K_{n-k}, (CK_{n-k})$. 
\end{definition}
Using the above definition, we prove the following theorem:
\begin{theorem}
\label{criterion for k prolific in Kn}
A permutation $\sigma$ is $k$-prolific in $K_n$ if and only if $br(\sigma)\geq k+3$.
\end{theorem}
\begin{proof}
Assume to the contrary that $\sigma$ is $k$-prolific but $br(\sigma)< k+3$. This means that there are $1\leq i,j \leq n$ such that $d_{\sigma}(i,j)<k+3$. For the sake of convenience, we assume that $i<j$ and $\sigma_i<\sigma_j$. 

Now, let \[S=\{\sigma_{m}|\,i<m<j \text{ or }  \sigma_i<\sigma_m<\sigma_j\}.\] Obviously, $|S|\leq j-i-1+\sigma_j-\sigma_i-1$ so that $|S|+2\leq d_{\sigma}(i,j)< k+3$ and thus $|S|\leq k$. 
Omitting the elements of $S$ from $\sigma$, we get a permutation $\sigma' \in S_{n-|s|}$ such that $\sigma' \prec \sigma $ and $\sigma'$ contains a pair $u<v$ such that $d_{\sigma'}(u,v)=2$, i.e. $(\sigma_u,\sigma_v)$ is a bond. 
For example, let $\sigma=[{\bf 9}72{\bf 4}15836]$ and let $i=1$ and $j=4$. Then we have $S=\{7,2,5,6,8\}$ and after omitting the elements of $S$, we get $\sigma'=[{\bf 43}12]$ which indeed has a bond.

Now, if $|S|=k$, this means that $\sigma'\notin K
_{n-k}$ so that $\sigma$ can not be $k-$prolific in $K_n$. 
If $|S|<k$ then we can omit from $\sigma'$ a set of $k-|S|$ elements which does not contain the bond mentioned above and  get an $(n-k)$-subset of the entries of $\sigma$ which creates a descendant of $\sigma$ which is not a member of $K_{n-k}$, so again $\sigma$ can not be $k-$prolific.\\

On the other hand, assume that $br(\sigma)\geq k+3$ and we prove that $\sigma$ is $k$-prolific in $K_n$. 
First, by Theorem 2.25 of \cite{BHT}, $br(\sigma)\geq k+2$ if and only if $\sigma$ is $k-$ prolific in $S_n$, so that each omission of a $k-$ set from $\sigma$ yields a  different permutation in $S_{n-k}$. 
To show that each such $\pi\in S_{n-k}$ is indeed a king permutation, note that  
by Proposition 2.24 of \cite{BHT}, omitting a single entry from a permutation decreases the breadth by at most one. Thus, omitting each $k-$ set of elements forms a permutation $\pi$ with $br(\pi)\geq 3$ which implies that $\pi\in K_{n-k}$. 


\end{proof}


The following theorem is the generalization of Theorem \ref{criterion for k prolific in Kn} for the poset of cylindrical king permutations. 

\begin{theorem}\label{criterion for k prolific in CKn}
A permutation $\sigma$ is $k$-prolific in $CK_n$ if and only if  $cbr(\sigma)\geq k+3$.
\end{theorem}

\begin{proof}[Sketch of proof]
The proof is almost identical to the proof of Theorem \ref{criterion for k prolific in Kn} except for the following necessary changes:  

\begin{itemize}
    \item  Instead of the set $S$, we define the set:
    \[S_{cyclic}=\{\sigma_{m}|\,m \in ((i,j)) \text{ or }  \sigma_i<\sigma_m<\sigma_j\}.\]

\item Replace the breadth $br$ by its cyclic  parallel
breadth $cbr$.

\item At the end of the proof, use Proposition \ref{single entry} of this article instead of Proposition 2.24 of \cite{BHT}.

\end{itemize}

\end{proof}


\section{The structure of the poset of cylindrical king permutation}
\label{structure}

It is easy to observe that for $\sigma \in S_n$, and a regular bond $(a,a+1)$ in $\sigma$, omitting $a$ or $a+1$ amounts to the same permutation. In other words  \begin{equation}\label{omit from bond}
   \nabla_a^*(\sigma)=\nabla_{a+1}^*(\sigma).
\end{equation}
\begin{example}
For $\sigma=[523641]$, $\nabla_2^*(\sigma)=\nabla_3^*(\sigma)=[42531].$
\end{example}
 This enables one to count for each permutation $\sigma \in S_n$, the number of permutations $\pi \in S_{n-1}$ such that $\pi \prec \sigma$.
 


Unfortunately, Equation (\ref{omit from bond}) does not hold anymore when we deal with edge bonds. For example, if $\sigma=[246351]$, then $\nabla_2^*(\sigma)=[35241]$ while $\nabla_1^*(\sigma)=[13524]$.
In order to overcome this problem, 
we define an equivalence relation based on the view of  each permutation as if it was written on a circle and not on a line. For example, according to this new view, the permutations $[35241]$ and $[13524]$ are considered identical. Formally, we define: 

\begin{definition}
Let $C_n$ be the cyclic sub-group of $S_n$, generated by the cycle $(12\cdots n)$.  Define an equivalence relation on $S_n$ by $\sigma \sim \tau$ if and only if there is some $\pi\in C_n$ such that $\tau =\sigma \cdot \pi$ and let $S_n/C_n$ be the quotient space. Also, let us denote by $\orb(\sigma)$ the equivalence class of $\sigma$. 
\end{definition}

\begin{example}
Let $\sigma=[23154] \in S_5$. Then $$\orb(\sigma)=\{[23154],[31542],[15423],[54231],[42315]\}.$$
\end{example}

The following theorem characterizes the permutations that serve as building blocks for the poset of cylindrical king permutations.
Note that, if $\sigma \sim \tau$ then $\sigma$ and  $\tau$ have the same number of cyclic bonds. Therefore, $\sigma \in CK_n$ iff $\tau \in CK_n$.

\begin{theorem}\label{Each king contains 31425 or 24135}
For every $\sigma  \in CK_n$ ($n \geq 5$) there is some $\pi \in \orb([31425])\cup \orb([24135])$ such that $\pi \preceq \sigma$.
\end{theorem}

\begin{proof}
We prove by induction on $n$, the case  $n=5$ being trivial, since $CK_5=\orb([31425])\cup \orb([24135])$. 
We assume that each cylindrical king permutation of order $n-1$ contains at least one of the permutations of $\orb([31425])\cup \orb([24135])$ and prove that each  cylindrical king permutation of order $n$ contains at least one element of this set.

Let $n > 5$ and assume to the contrary that there is $\sigma=[\sigma_1,\dots,\sigma_n] \in CK_n$ which does not contain any of the permutations in $\orb([31425])\cup \orb([24135])$. Then for every $1 \leq i \leq n$,  the permutation $\nabla^*_i(\sigma)$,  contains neither of them as well. By the induction hypothesis, $\nabla^*_i(\sigma) \notin CK_{n-1}$, which implies, that every digit of $\sigma$ is a cyclic separator. In particular, $1$ is a cyclic vertical separator.
Hence, we have: $\sigma\in \orb([\ldots, a,1,a+1,\ldots])$ or $\sigma\in \orb([\ldots, a+1,1,a,\ldots])$. Since $\sigma \in CK_n$, we know that  $a>2$.
Without loss of generality, assume that  $\sigma\in \orb([\ldots, a,1,a+1,\ldots])$. 
By the above, the element $a+1$ is a cyclic separator of $\sigma$. We have now two options as $a+1$ is horizontal or vertical :
\begin{itemize}
    \item If $a+1$ is  horizontal then $\sigma \in \orb([\ldots a+2,a,1,a+1,\ldots])$. For the same reasons as above, $a$ is also a
        cyclic horizontal separator of $\sigma$.  Thus  $\sigma \in \orb([\ldots a+2,a,1,a+1,a-1\ldots])$. But the sub sequence $a+2,a,1,a+1,a-1$ is a pattern which is isomorphic to $[53142]$, so that $\sigma$ contains an element of $\orb([31425])$.
        
        \item If $a+1$ is vertical then $\sigma\in \orb([\ldots a,1,a+1,2\ldots])$. As we saw above, $a$ is also a cyclic separator of $\sigma$ which must be horizontal so that $a=3$. This implies that  $\sigma\in \orb([\ldots 3,1,4,2\ldots])$. Since the element $5$ has to be somewhere in $\sigma$ thus  $\sigma \in \orb([\ldots ,3,1,4,2,\ldots,5,\ldots])$ which contains the pattern $[31425]$.
\end{itemize}

\end{proof}

In proposition 3.14 of \cite{BERSKINGS}, it is proven that for each two king permutations $\sigma$ and  $\pi$ such that $\pi \prec \sigma$ and $|\sigma|-|\pi|=2$ there is a king permutation $\tau$ such that $\pi \prec \tau \prec \sigma$. This does not hold for cylindrical king permutations. For example, for $\sigma=[18463527]\in CK_8$ and $\pi=[635241]\in CK_6$, even though $\pi \prec \sigma$, there is no $\tau \in CK_7$ such that $\pi \prec \tau \prec \sigma$.

Moreover, in Theorem 3.15 of \cite{BERSKINGS} it is claimed  that if $\sigma,\pi$ are king permutations such that $\pi \prec \sigma$ and $|\sigma|-|\pi|>3$ then there exists a king permutation $\tau$ such that $\pi \prec \tau \prec \sigma$ and $|\sigma| -|\tau| \in \{1,3\}$.
I.e. the maximal gap between any two king permutations in the poset $\mathbb{K}$ is $3$. This property does not hold true anymore in the poset $\mathbb{CK}$ as the following example verifies:

\begin{example}

Let $\sigma=[579683142]$ and $\pi=[13524]$ which is obtained from $\sigma$  by removing the elements $3,1,4,2$. It is easy to observe that there is no $\tau \in \mathbb{CK}$ such that $\pi \prec \tau \prec \sigma$.

\end{example}

We claim now that the maximal gap between two permutations in $\mathbb{CK}$ is $4$. 

\begin{theorem}\label{antisematic}

Let $\sigma, \pi \in \mathbb{CK}$ be such that $\pi \prec \sigma$ and $|\sigma|-|\pi|>4$.
Then there is some $\tau \in \mathbb{CK}$ such that $\pi \prec \tau \prec \sigma$ and $|\sigma|-|\tau|\leq 4$.

\end{theorem}
    
The proof of Theorem \ref{antisematic} is very subtle and messy. It consists of many similar sub-cases which recur over and over. In order not  to be too exhaustive, we present some basic observations and lemmas and then show how to deal with few typical sub-cases and leave the rest of the proof for the reader.   
For the sake of simplifying the proof, we present a series of claims:

\begin{observation}\label{Lemma 1}
Let $\sigma$ be a permutation, let $\pi \in \mathbb{CK}$, and let $i$ be such that $\pi \prec \nabla^*_i(\sigma)$. If $\nabla^*_i(\sigma)$ contains the cyclic bond consisting of the elements $j$ and $k$ (of $\sigma$) then $\pi \prec \nabla^*_k(\sigma)$ or $\pi \prec \nabla^*_j(\sigma)$, 

Moreover, if the bond is not an edge bond then both $\pi \prec \nabla^*_k(\sigma)$ and $\pi \prec \nabla^*_j(\sigma)$.
\end{observation}
\begin{example}
In the left picture of Figure \ref{S2}, $\sigma=[264159837]\in S_9, \pi=[241635]\in CK_6$, such that $\pi \prec \sigma$,
$\pi \prec \nabla^*_5(\sigma)$.
In $\nabla^*_5(\sigma)=[25418736]$ there is  a regular bond  $(5,4)$ (the digits $\textcolor{purple}{4}$ and $\textcolor{purple}{6}$ in $\sigma$),  we get that $\pi \prec \nabla^*_4(\sigma)$ and $\pi \prec \nabla^*_6(\sigma)$.

In the right picture $\sigma=[264159837]\in S_9, \pi=[253146]\in CK_6$, such that
$\pi \prec \sigma$,
$\pi \prec \nabla^*_7(\sigma)$.
In $\nabla^*_7(\sigma)=[26415873]$ there is an edge bond $(2,3)$ (the digits $\textcolor{purple}{2}$ and $\textcolor{purple}{3}$ in $\sigma$) ,  we get that $\pi \prec \nabla^*_3(\sigma)$ but $\pi \nprec \nabla^*_2(\sigma)$ .
\end{example}

\begin{figure}[!ht]
    \centering
\begin{tikzpicture}
        \tikzstyle{every node} = [rectangle]
        
       \node (264159837) at (0,0) {$\textcolor{purple}{\sigma=[264159837]}$};
        
        \node at (0.6,-1) {$5$, ($\textcolor{purple}{5}$)};
        \node (25418736) at (0,-2) {$[25418736]$};
        
        \node at (1.1,-3) {$4$ or $5$ ($\textcolor{purple}{4/6}$) };
        \node (2417635) at (0,-4) {$[2417635]$};
        
        \node at (1.1,-5) {$6$ or $7$ (\textcolor{purple}{8/9})};
        \node (241635) at (0,-6) {$\pi=[241635]$};

        
        \node (264159837b) at (8,0) {$\textcolor{purple}{\sigma=[264159837]}$};
        
        \node at (8.6,-1) {$7$ ($\textcolor{purple}{7}$)};
        \node (26415873) at (8,-2) {$[26415873]$};
        
        \node at (9.1,-3) {$7$ or $8$ ($\textcolor{purple}{8/9}$) };
        \node (2641573) at (8,-4) {$[2641573]$};
        
        \node at (9.3,-5) {$3$ (\textcolor{purple}{3})};
        \node (253146) at (9.3,-6) {$\pi=[253146]$};
        
        \node at (6.8,-5) {$2$ (\textcolor{purple}{2})};
        \node (531462) at (7,-6) {$[531462]$};

        \foreach \from/\to in {264159837/25418736, 25418736/2417635,2417635/241635,26415873/2641573,2641573/253146,
        2641573/531462,264159837b/26415873} 
            \draw[->] (\from) -- (\to);
    \end{tikzpicture}
    \caption{}
    \label{S2}
\end{figure}

\begin{definition}
Let us denote by $\nabla^*_A(\sigma)$ the permutation $\nabla^*_{i_1}(\nabla^*_{i_2}( \dots(\nabla^*_{i_n}(\sigma)))$ where $A=\{i_1,\dots,i_n\}$ and $i_1<\dots<i_n$.
For example $\sigma=[264159837]\in S_9, \pi=[241635]\in CK_6$,
 so that $\pi \prec \sigma$.
 For $A_1=\{4,5,8\},A_2=\{4,5,9\}, A_3=\{5,6,8\}, A_4=\{5,6,9\}$, we get $\pi=\nabla^*_{A_i}(\sigma)$ for $1 \leq i \leq 4$.
 \end{definition}



In the following observation, when we write 'adjacent', we take it in the cyclic meaning, i.e. in $[a_1,\dots,a_n]$, the elements $a_1$ and $a_n$ are considered here as adjacent.

\begin{observation}
\label{obs}

     Let $\sigma \in \mathbb{CK}$ such that $\sigma$ contains the (cyclic) block $(a+2,a+4,a+1,a+3)$ or its reverse, i.e. $\sigma \in orb([\dots,a+2,a+4,a+1,a+3,\dots])$ or $\sigma \in orb([\dots,a+3,a+1,a+4,a+2,\dots])$. 
     \begin{enumerate}
     \item If we remove any single element $i$ from the block $(a+2,a+4,a+1,a+3)$ in $\sigma$ then $\nabla^*_i(\sigma)$ will contain a cyclic bond. Moreover, we will not obtain a cylindrical king permutation out of $\sigma$ unless we omit at least $3$ elements from the block, in which case, it contracts to the element $a+1$.   
    
     \item  If the elements $a$ and $a+5$ are not placed adjacent to the block $(a+2,a+4,a+1,a+3)$ or its reverse in $\sigma$, then by omitting any $3$ elements of the block we get a cylindrical king permutation.    
    
     For example, let $\sigma=[39142\bf{6857}]$ with the block $(6,8,5,7)$. Note that $\tau=\nabla_6^*(\sigma)=[38142756]$ which contains the bond $(5,6)$, $\nabla_7^*(\tau)=[3714256]$ with the bond $(5,6)$, but $\nabla^*_{\{5,6,7\}}(\sigma)=[361425]$ which is a cylindrical king permutation. 
     
    \item If $a$ or $a+5$ is placed adjacent to the block $(a+2,a+4,a+1,a+3)$ or its reverse in $\sigma$, then removing it produces a cylindrical king permutation. This is due to the fact that in this case $a$ and $a+5$ can not be cyclic separators of any type in $\sigma$.   
    For example, let $\sigma=[3192{\color{red}{4}} {\bf{6857}}  ]$ with the block $(6,8,5,7)$. Note that $\tau=\nabla_4^*(\sigma)=[31825746]\in CK_{8}$. 
    
\end{enumerate}
\end{observation}

\begin{lemma}\label{second lemma}
Let $\sigma,\pi \in \mathbb{CK}$ be such that $\pi \prec \sigma$ and $|\sigma|-|\pi|>4$. Assume that 
\begin{itemize}
    \item 
    There is some element $a\in \mathbb{N}$ such that $\sigma \in orb([\ldots,a+2,a+4,a+1,a+3,\ldots])$ or $\sigma \in orb([\ldots,a+3,a+1,a+4,a+2,\ldots])$ . 
    \item There is some $1 \leq i \leq 4$ such that $\pi \prec \nabla^*_{a+i}(\sigma)$.
     
\end{itemize}

Then there is some $\tau \in \mathbb{CK}$ such that $\pi \prec \tau \prec \sigma$
and $|\sigma|-|\tau|\leq 4$.  
\end{lemma}

\begin{proof}
Assume without loss of generality that  
$\pi \prec \nabla^*_{a+2}(\sigma)$ and that $\sigma$ contains the block $(a+2,a+4,a+1,a+3)$. 
The proof depends on the location of this block in $\sigma$. The possibilities are: 
$\sigma=[\ldots ,a+2,a+4,a+1,a+3]$ or $\sigma=[a+2,a+4,a+1,a+3,\ldots]$ or $\sigma=[\ldots,a+2,a+4,a+1,a+3,\ldots]$ or $\sigma=[a+1,a+3,\ldots,a+2,a+4]$ or $\sigma=[a+3,\ldots,a+2,a+4,a+1]$ or $\sigma=[a+4,a+1,a+3,\ldots,a+2]$. 

Here we show only the first case: $\sigma=[\dots ,a+2,a+4,a+1,a+3]$ which seems to be the most general case in the sense that it contains all the arguments we need. The other cases can be proved using the same method. 
    

According to Observations \ref{Lemma 1} and \ref{obs}$.1$, $\pi \prec \nabla^*_A(\sigma)$ for $A=\{a+2,a+3,a+4\}$.
If $\nabla^*_A(\sigma) \in \mathbb{CK}$ then we are done (just take $\tau=\nabla^*_A(\sigma)=[\dots,a+1]$). 
 Otherwise,  according to \ref{obs}.2, either $a$ or $a+5$ is adjacent to the the block ($a+2,a+4,a+1,a+3$). Assume w.l.o.g. that $a$ is adjacent to the above block. We have $2$ possibilities. 
 \begin{itemize}
     \item  If $\sigma=[\dots,a,a+2,a+4,a+1,a+3]$ then $\nabla^*_A(\sigma)=[\dots,a,a+1]$ and according to Observation \ref{Lemma 1}, $\pi \prec \nabla^*_a(\nabla^*_A(\sigma)) = \nabla^*_{A'}(\nabla^*_a(\sigma)) \prec \nabla^*_a(\sigma)$ for $A'=\{a+1,a+2,a+3\}$ and so by \ref{obs}.3 $\tau=\nabla^*_a(\sigma)$ and we are done.
     \item If $\sigma=[a,\dots,b,a+2,a+4,a+1,a+3]$ (where $b \neq a+5$) then $\nabla^*_A(\sigma)=[a,\dots,a+1]$. 
 By Observation \ref{Lemma 1}, $\pi\prec \nabla_a^*(\nabla^*_A(\sigma))\prec \nabla_a^*(\sigma)$ or  $\pi \prec \nabla^*_{B}(\sigma)$ for $B=A\cup{\{a+1\}}$. 
\begin{enumerate}
    \item  If $\pi\prec \nabla^*_a(\sigma)$ then we are done since according to \ref{obs}.$3.$,  $\tau=\nabla_a^*(\sigma)$ is a cylindrical king permutation. 
 \item If  $\pi\nprec \nabla_a^*(\sigma)$ and $\nabla_B^*(\sigma)$ is a cylindrical king permutation then we are done by taking  $\tau=\nabla^*_B(\sigma)$. (Note that this is the single case in which $|\sigma|-|\tau|=4$).
 Otherwise, the permutation $\nabla^*_B(\sigma)$ contains a cyclic bond, leading to one of the following two possibilities:
 \begin{enumerate}
     \item $\sigma=[a,a+5,\ldots,a+2,a+4,a+1,a+3]$ which is impossible since then $\nabla^*_B(\sigma)$ is of the form $[a,a+1,\dots]$ and by \ref{Lemma 1} we get $\pi\prec\nabla^*_{a}(\sigma)$ which contradicts the assumption.
     \item $b=a-1$, i.e. $\sigma=[a,\ldots,a-1,a+3,a+1,a+4,a+2]$ and according to Observation \ref{Lemma 1}, $\pi\prec\nabla_{a-1}^*(\nabla^*_A(\sigma))\prec\nabla_{a-1}^*(\sigma)$ . We claim that $\nabla_{a-1}^*(\sigma)$ is a cylindrical king permutation. Indeed, otherwise $\sigma=[a,a-2\dots,a-1,a+2,a+4,a+1,a+3]$ and $\nabla_{a-1}^*(\sigma)$ contains a regular bond consisting of the elements $a$ and $a-2$ of $\sigma$. Hence, due to \ref{Lemma 1} $\pi\prec \nabla^*_{a}(\sigma)$ which contradicts the assumption.
Now, we can take $\tau=\nabla_{a-1}^*(\sigma)$ and we are done.  

 \end{enumerate}  
\end{enumerate}
\end{itemize}
 

\end{proof}
Now we return to the proof of Theorem \ref{antisematic}. Let $\sigma, \pi$ be cylindrical king permutations such that $\pi \prec \sigma$, and $|\sigma|-|\pi|>4$ and let us prove that there is $\tau \in \mathbb{CK}$ such that $\pi \prec \tau \prec \sigma$ and $|\sigma|-|\tau|\leq 4$.
\begin{proof}{of Theorem 3.6}

Since $\mathbb{CK}\subset \mathbb{K}$, we have that $\pi$ and $\sigma$ are king permutations. Thus, according to Theorem 3.15 in \cite{BERSKINGS},
there is $\tau \in \mathbb{K}$ such that $\pi \prec \tau \prec \sigma$ and $|\sigma|-|\tau|\in\{1,3\}$. 
If $|\sigma|-|\tau|=3$ then according to the proof of Theorem 3.15 in \cite{BERSKINGS}, $\sigma$ contains the block $(a+2,a+4,a+1,a+3)$ or its reverse and $\tau$ is obtained from $\sigma$ by omitting $3$ elements from this block. Then we can proceed according to Lemma \ref{second lemma}. 
Assume now that $|\sigma|-|\tau|=1$. If $\tau \in \mathbb{CK}$ then we are done. Otherwise, $\tau$ is of the form: $\tau=[a,\dots,a+1]$ and there exists some cyclic separator $b$ such that $\tau=\nabla^*_b(\sigma)$. If $b$ is vertical then we have $\sigma=[a,\dots,a+1,b]$, otherwise 
$b$ is horizontal i.e. $b=a+1$ and so
$\sigma=[a,\dots,a+1,\dots,a+2]$ (or their reverses). Here we consider only the vertical case: $\sigma=[a,\dots,a+1,b]$. 
Since $\pi \prec \nabla^*_b(\sigma)=[a,\dots,a+1]$, by Observation \ref{Lemma 1}, $\pi \prec \nabla^*_
a(\sigma)$ or $\pi \prec \nabla^*_{a+1}(\sigma)$. 

We consider here only the case $\pi \prec \nabla^*_{a+1}(\sigma)$ but $\pi \nprec \nabla^*_a(\sigma)$. If $\nabla^*_{a+1}(\sigma)\in \mathbb{CK}$ we are done. Otherwise  $\nabla^*_{a+1}(\sigma)$ contains a cyclic bond. Thus we have for $\sigma$ one of the following options: 
either $\sigma=[a,\dots,b\pm 1,a+1,b]$ or $\sigma=[a,a+2\dots,a+1,b]$. Let us consider again the first case, $\sigma=[a,\dots,b+1,a+1,b]$. In this case, $\nabla^*_{a+1}(\sigma)$ has a regular bond consisting of the elements $b$ and $b+1$ of $\sigma$ and so by Observation \ref{Lemma 1}, $\pi \prec \nabla^*_b(\sigma)$ and $\pi \prec \nabla^*_{b+1}(\sigma)$. 
If $ \nabla^*_{b+1}(\sigma)$ is a cylindrical king permutation, we are done. Otherwise, the omission of the element $b+1$ produces a cyclic bond. Again, there are two options for the structure of $\sigma$. The first option is that $a=b+2$ and thus $\sigma=[a,\dots,a-1,a+1,a-2]$ and thus according to Lemma \ref{second lemma} we are done.
The second option (if $a\neq b+2$) is that $\sigma=[a,\dots,a+2,b+1,a+1,b]$, in this case $\nabla^*_A(\sigma)=[a,\dots,a+1,b]$ for $A=\{a+2,b+1\}$ is a cylindrical king permutation. (Indeed, the pair $(a+1,b)$ is not a bond since $a+1$ and $b$ are adjacent in $\sigma$. The elements $a+1$ or $a-1$ can not appear beside the element $a$ since $\sigma$ is a cylindrical king permutation). This completes the proof. 
\end{proof}

\section{Directions for further research}
\label{direction} 
\begin{itemize}
    \item In \cite{BHT},  Bevan,  Homberger and Tenner proved that a $k$-prolific permutation in $S_n$ must be at least of size $\frac{k^2}{2}+2k+1$. 
A challenge may be to find a lower bound for the size of a $k$-prolific permutation in $K_n$ and in $CK_n$.
\item In the chess game, one can further enable the pieces to go off a column and reappear at the beginning of that columns, thus obtaining a 'chess on the torus'. The permutations describing non-attacking kings on this board are called torical king permutations and constitute a sub-poset of their own. The work on these permutations is in progress. 
\end{itemize}

\section{Appendix}

In a recent paper \cite{separator}, the authors defined a new concept, called {\it separator}. A separator in a permutation $\sigma=[\sigma_1,\dots,\sigma_n] \in S_n$ is an element $\sigma_i$
that its omission forms a new bond. In particular, if $\sigma \in K_n$ is a king permutation, then the permutation obtained by the omission of $\sigma_i$ is a king permutation if and only if $\sigma_i$ is not a separator in $\sigma$. In order to characterize the structure of the poset of cylindrical kings, we need a cyclic version of that concept. We start by recalling the definition of the (regular) separator (see Definition 2.1 in \cite{separator}). 
\begin{definition}\label{def separate}
For $\sigma=[\sigma_1,\ldots,\sigma_n] \in S_n$ we say that $\sigma_i={\bf\it \textcolor{red}{a}}$ {\em separates} $\sigma_{j_1}$ from $\sigma_{j_2}$ in $\sigma$ if by omitting the element $a$ from $\sigma$ we get a {\bf new} bond. This happens if and only if one of the following cases holds: 

\begin{enumerate}
\item $j_1,i,j_2$ are subsequent numbers and $|\sigma_{j_1}-\sigma_{j_2}|=1$, i.e $$\sigma=[\ldots,{\bf \it b},{\bf \it \textcolor{red}{a=\sigma_i}},{\bf \it b \pm 1},\ldots]$$
In this case we say that $a$ is a {\bf vertical separator}.
\item $\sigma_{j_1},\sigma_i=a,\sigma_{j_2}$ is an increasing or decreasing sub sequence of subsequent numbers, and $|j_1-j_2|=1$, i.e, $$\sigma=[\ldots,{\bf \it \textcolor{red}{a}},\ldots,{\bf \it a \pm1,a\mp1},\ldots]$$ or $$\sigma=[\ldots,{\bf \it a\pm1,a\mp1},\ldots,{\bf \it \textcolor{red}{a}},\ldots].$$
In this case we say that $a$ is a {\bf horizontal separator}.
\end{enumerate}
\end{definition}

Cylindrical kings are actually permutations having no {\bf cyclic} bonds. In order to be able to deal with edge bonds, we have to modify the definition \ref{def separate} of the separator.

\begin{definition}
 Let $\sigma=[\sigma_1,\ldots,\sigma_n] \in S_n$. An element $\sigma_i$ will be called an {\bf edge separator} if omitting it produces a new edge bond. This happens when one of the following holds: 
\begin{enumerate}
\item  If $|\sigma_1-\sigma_{n-1}|=1$ then $\sigma_n$ is a (vertical) edge separator. 
\item If $|\sigma_n-\sigma_2|=1$ then $\sigma_1$ is a (vertical) edge separator.
\item If for some $1<i<n$ we have that $\sigma_1,\sigma_i,\sigma_n$ is an increasing or decreasing sub sequence of subsequent numbers, then we say that $\sigma_i$ is a (horizontal) edge separator.
\end{enumerate}
We say that an element $\sigma_i$ is a {\bf cyclic separator} if it is a regular or an edge separator. In this case we say that $\sigma_i$ cyclically separates some element from another.
\end{definition}

\begin{example}
Let $\sigma=[52341]$, then $\sigma_5=1$ is a (vertical) edge separator in $\sigma$. Let $\sigma=[6257134]$, then $\sigma_3=5$ is a (horizontal) edge separator in $\sigma$. 
If $\sigma=[246135]$, then $\sigma_6=5$ is  a vertical edge separator and a horizontal (regular) separator in $\sigma$.
\end{example}

\printbibliography

\end{document}